\documentclass[12pt]{amsart}
\usepackage{amscd,amssymb}

\numberwithin{equation}{section}

\newtheorem {Theorem} 			{Theorem}
\newtheorem {varTheorem}                {Theorem}
\newenvironment {Theorem'}
        {\begin{varTheorem}{\hspace{-3.5mm}}{\bf '}{\hspace{3.5mm}}}
        {\end{varTheorem}}
\newtheorem {RefTheorem}[equation]     	{Theorem}  	
\newtheorem {Lemma}[equation]     	{Lemma}  	
\newtheorem {Proposition}[equation]	{Proposition}

\theoremstyle{definition}

\newtheorem {Remark}[equation]		{Remark}

\newcommand{\pr} {\smallskip\noindent{\bf Proof\,\,}}

\newenvironment{pf}  {\begin{proof}}{\end{proof}}

\newcommand     {\comment}[1]   {}
\newcommand     {\mute}[2] {}
\newcommand     {\printname}[1] {}

\newcommand{\labell}[1] {\label{#1}\printname{#1}}

\def	\red	{{\operatorname{red}}}

\def    \ker	{\operatorname{ker}}

\def	\inv	{^{-1}}
\def	\to	{\longrightarrow}

\def	\ssminus	{\smallsetminus}

\def	\C	{{\Bbb C}}
\def	\R	{{\Bbb R}}
\def	\bN	{{\Bbb N}}
\def	\Q	{{\Bbb Q}}
\def	\Z	{{\Bbb Z}}
\def	\CP	{{\Bbb C}{\Bbb P}}
\def	\RP	{{\Bbb R}{\Bbb P}}
\def	\Mred	{{M_\red}}

\def	\cF	{{\mathcal F}}

\def	\cI	{{\mathcal I}}
\def	\cJ	{{\mathcal J}}

\def	\t	{{\frak t}}
\def	\ft	{{\frak t}}
\def	\fg	{{\frak g}}

\def	\pr	{\operatorname{pr}}
\def	\ssminus	{\smallsetminus}

\begin{document}


\title[The cohomology rings of abelian symplectic quotients]
{The cohomology rings of abelian symplectic quotients}

\author{Susan Tolman}
\author{Jonathan Weitsman}
\thanks{S. Tolman was partially supported by an NSF Mathematical Sciences
Postdoctoral Research Fellowship.
J. Weitsman was supported in part by NSF grant DMS 94/03567, by NSF 
Young Investigator grant DMS 94/57821, and
by an Alfred P. Sloan Foundation Fellowship.}

\address{Department of Mathematics,  University of Illinois at
Urbana-Champaign, 
Urbana, IL 61801}
\email{stolman@math.uiuc.edu}

\address{Department of Mathematics, University of California, Santa Cruz,
CA 95064}
\email{weitsman@cats.ucsc.edu}
\thanks{\today}

\begin{abstract}
Let $M$ be a symplectic manifold, equipped with
a Hamiltonian action of a torus $T$.  We give an explicit formula
for the rational cohomology ring of the symplectic quotient $M//T$
in terms of the cohomology ring of $M$ and fixed point data.
Under some restrictions, our formulas apply to integral 
cohomology.  In certain cases these methods enable us to show
that the cohomology of the reduced space is torsion-free. 
\end{abstract}

\maketitle



\section{Introduction}

Let a torus $T$ act on  a compact symplectic manifold $(M,\omega)$
with a moment map $\phi:M\to {\frak t}^* = {\rm Lie}(T)^*$.
Assume that $0$ is a regular value of $\phi$, and  let
$\Mred := \phi^{-1}(0)/T$  denote the reduced space.
The inclusion $i:\phi^{-1}(0)\to M$ induces a map in equivariant
cohomology, which, when composed with the natural identification
$H^*_T(\phi^{-1}(0)) \simeq H^*(\Mred)$, yields a map $\kappa:H^*_T(M)
\to H^*(\Mred)$, called the {\em Kirwan map}.  According to a theorem
of  Kirwan, this map is a surjection.  The purpose of this
paper is to answer the natural question:  What is the kernel
of $\kappa$?

This question was answered explicitly by Kirwan in a number of
examples \cite{K1}.
In \cite{witten} Witten formulated the idea of nonabelian localization: this
term refers to a procedure which computes the evaluation on the
fundamental class of a symplectic quotient 
$\int_{M_{\rm red} } \kappa (\alpha) $ (for any equivariant
cohomology class $\alpha$) in terms of data on the original 
manifold $M$. For Witten the data on $M$ consisted of the
critical sets of the function $||  \phi ||^2$ (the norm
squared of the moment map). Motivated by Witten's work,
Jeffrey and Kirwan [JK] gave a proof of  Witten's results using the
geometry of the image of the moment map, and also found an alternative
version of nonabelian localization (the residue formula,
Section 8 of [JK]) which computes 
$\int_{M_{\rm red} } \kappa (\alpha) $ in terms of 
fixed point data on $M$ (in other words the components of the
fixed point set of the maximal torus, the characteristic classes 
of their normal bundles, the weights of the action of the
maximal torus on the normal
bundles, and the values of the moment map on the components
of the fixed point set).
Since, by Poincar\'e duality, 
$\beta \in H_T^*(M;\Q)$ is in the kernel of $\kappa$
exactly if the integral of $\kappa( \alpha \beta)$ 
over $M_\red$ is  zero for all $\alpha \in H_T^*(M,\Q)$, 
the kernel of $\kappa$ can, in principle, be computed from
their results, in the case of rational cohomology.
However, we do not know of any direct method of relating
their formulas to ours.  See also \cite{gk,p,v}.

In this paper we give a description
of the kernel of $\kappa$ in terms of fixed point data.
For according to a theorem of F. Kirwan, the natural restriction
map from the equivariant cohomology of $M$ to the equivariant cohomology
of the fixed point set is  an injection.
We compute the kernel of $\kappa$ in
terms of the image of this map.
This image is well understood in many examples.
Moreover, it is determined by the subset of one and zero dimensional orbits
(see \cite{GKM, TW} and references therein).
Our methods also give a basis for the kernel of $\kappa,$
and, under some restrictions, allow us to compute
the {\em integral} cohomology
rings of symplectic quotients.

Our main result is the following.

{\bf Theorem 2.} {\em  Let a torus $T$ act on a compact symplectic
manifold $M$ with a
moment map $\phi: M \to \t^*$.  Assume that $0$ is a regular value of $\phi$.
Let $F$ denote the set of fixed points.
For all $\xi \in \ft$, define
$$M_\xi := \{ m \in M \mid  \left< \phi(m), \xi \right> \leq 0 \} ,$$
$$K_\xi := \{ \alpha \in H^*(M;\Q) \mid \alpha|_{F \cap M_\xi} = 0\  \},
\hbox{and} $$ 
$$K := \sum_{\xi \in \ft} K_\xi.$$
Then there is a short exact sequence:
$$ 0 \to K \to H_T^*(M;\Q) \stackrel{\kappa}{\to}
 H^*(M_\red;\Q) \to 0,$$
where 
$\kappa: H_T^*(M;\Q) \to H^*(\Mred;\Q)$ is  the Kirwan map.
}

The main ingredient in the proof of this theorem is an application
of Morse theory (as extended by Bott and Kirwan) to functions associated
to the moment map $\phi$.  Specifically, we study the functions given
by the components $\phi^\xi$ of the moment map, and by its square
$||\phi||^2$.  The Morse theory of these functions has remarkable
properties, as shown by Frankel (\cite{F}) and by
Atiyah and Bott \cite{AB,AB2}.  These properties can be summarized
by the statement that as Morse functions, these functions are
self-completing, and equivariantly, they are locally self-completing.
More technically, we choose some critical value $c$ of the relevant function
$f$, and suppose that the interval $[c-\epsilon, c+\epsilon]$ contains
no critical values of $f$ other than $c$.  Write 
$M_\pm = f^{-1}(-\infty,c\pm \epsilon)$.  Then the long exact sequence
of the pair $(M_+, M_-)$ breaks up into short exact sequences

$$0\to H^*(M_+,M_-)\to H^*(M_+) \to H^*(M_-)\to 0 \label{exactness},$$

\noindent where the $H^*(\cdot)$ denotes either ordinary or equivariant
rational cohomology.  

When applied to the function $f = \phi^\xi$ for some $\xi$ in the
Lie algebra ot $T$, this sequence can be used to study the equivariant
cohomology of $T$.  For example, Kirwan used this
approach to prove that $f$ is a perfect and equivariantly
perfect Morse function, and that the cohomology of $M$ is equivariantly
formal.  More importantly for our purposes, Kirwan also
showed that the restriction of the equivariant cohomology
of $M$ to the equivariant cohomology of the fixed point set
in an injection.  
This technique is used in a companion paper to analyze the image
of this map \cite{TW}.

In order to study the symplectic quotient $\Mred$, we consider instead
the Morse theory associated to the function $f=||\phi||^2$.  This
is convenient because the minimum of this function is precisely
$\phi^{-1}(0)$, whose equivariant cohomology coincides with the (ordinary)
cohomology
of $\Mred$.  In the context of gauge theory such Morse functions were
studied
by Atiyah and Bott \cite{AB2}; in the 
context of finite-dimensional manifolds they were studied by
Kirwan \cite{K}.  This variant of Morse theory is
the key ingredient in the proof of the surjectivity
theorem cited above.  As it turns out it is also the key element in
our proof of Theorem 2.

The remainder of this paper is structured as follows.  In Section 2,
we study the key local lemma, due to Atiyah and Bott, which implies
the equivariant perfection of the various Morse functions which we
will study.  We then proceed, in Section 3, to give a simple proof
of a version (Theorem 1) of our main theorem which applies to the
case where the torus $T$ is one-dimensional.  This serves to illustrate
the basic ideas of our proof in a simple setting, where we may use
Morse-Bott theory rather than its elaboration to the case of degenerate
Morse functions.  We then bite the bullet in Section 4, introducing
the ideas of Morse-Kirwan theory that will be required in Section 5, where
the main theorem (Theorem 2) is proved; we include a version for noncompact
manifolds as Theorem 3.

The last three sections are devoted to extensions of the basic ideas of
the paper.  In Section 6, we state and prove versions of our results which
hold for cohomology with integer coefficients.  In Section 7 we use
methods similar
to those used in the proof of the main theorem to give a condition
for the cohomology rings of symplectic quotients to be torsion-free.
Finally in Section 8 we show how our method can be used to compute the
cohomology rings of compact smooth projective toric varieties
(which are all given as symplectic quotients by tori of Euclidean spaces).

{\bf Acknowledgement}:  We would like to thank Victor Guillemin, who encouraged
us to consider the consequences for symplectic geometry of the work of
Goresky, Kottwitz and MacPherson \cite{GKM}.  The present
paper and its companion
\cite{TW} would not have been possible without his prescient advice.
One of us (J.W.) would also like to thank Prof. Guillemin
and the MIT Mathematics
Department for the opportunity
to spend the Spring semester of 1996 at MIT.

\section{The key lemma}

The first step in the proof of our main theorem is a local lemma.
This lemma, due to Atiyah and Bott, is the key fact behind
many of the results we describe.

\begin{Lemma} {\bf (Atiyah-Bott)}
\label{Atiyah-Bott}
Let a compact Lie group  $G$ act on a manifold $X$
and an oriented real vector bundle $E$ over $X$.
Assume that a circle subgroup $S^1 \subset G$ acts on $E$
so that the fixed point set is precisely $X$.
Choose an invariant metric on $E$, and let $D$ and $S$ denote the unit disk and
sphere bundle, respectively.  
Then the natural long exact sequence in relative equivariant cohomology
splits into a short exact sequence:  
$$ 0 \to H_G^*(D,S;\Q) \to  H_G^*(D;\Q) \to H_G^*(S;\Q) \to 0.$$
\end{Lemma}

Since $D$ is homotopic to $X$, we have an identification of 
$H_G^*(D;\Q)$ and $H_G^*(X;\Q)$.
Moreover, the  Thom isomorphism
identifies  $H_G^{*-\lambda}(X;\Q) $ and $ H_G^{*}(D,S;\Q)$.
Under these identifications, the natural map   
$H_G^*(D,S;\Q) \to  H_G^*(D;\Q)$ in the long exact sequence in
relative cohomology 
is identified with the map from $H_G^{*-\lambda}(X;\Q)$ to
$H_G ^{*}(X;\Q)$
given by multiplication by the equivariant Euler class of $E$,
where $\lambda$ is the dimension of $E$.
An alternate formulation of this lemma is that this Euler class
is not a zero divisor.

Notice that the corresponding sequence in ordinary cohomology
is never exact for any vector bundle of positive rank for dimensional reasons.

\section{The simplest case: $S^1$ actions}

In this section we prove our theorem in a
special case. Specifically, we compute 
the rational cohomology ring of the
reduction of a compact symplectic manifold by a circle.
The essential features of our argument are already clear in this
case. 

\begin{Lemma}
\labell{KM-circle}
Let the circle $S^1$ act on a compact symplectic manifold $M$ with a
moment map $\phi : M \to \R$. 
Let a function $f: M \to \R$ be defined in any of the
three 
following ways:
let $f:= \phi^2$, let $f := \phi$,
or let $f := -\phi$.  

Let $C$ be a critical set of index $\lambda$ for $f$,
and assume  that
$C$ is the only critical set in $f^{-1}( f(C) - \epsilon, f(C) + \epsilon)$ 
for $\epsilon > 0$.
Define $M^+ := f\inv(-\infty,f(C) + \epsilon)$, and 
$M^- := f\inv(-\infty,f(C) - \epsilon)$.
Then there exists a short exact sequence in equivariant cohomology
$$ 0 \to H_{S^1}^{*-\lambda}(C;\Q) \to  H_{S^1}^{*} (M^+;\Q) 
\to H_{S^1}^{*}(M^-;\Q) \to 0.$$
Moreover, the composition of the injection 
$H_{S^1}^{* - \lambda}(C) \to H_{S^1}^{*}(M^+)$
with  the restriction map  
$H_{S^1}^{* }(M^+) \to H_{S^1}^{*}(C)$
is the cup product with the Euler class of the negative normal
bundle of $C$.
\end{Lemma}

In other words, every equivariant cohomology class on $M^-$ extends to
$M^+$. The restriction to $C$ of cohomology classes
on $M^+$ which vanish on $M^-$ is injective and its image is
the set of multiples of the Euler class of the negative normal
bundle of $C$.

\begin{pf}
Except in the case where $f = \phi^2$ and $C = \phi\inv(0)$, the
function $f$ is Morse-Bott at every connected component $C$ of
the critical set, the set 
$C$ is a component of the fixed point set, 
and the negative normal bundle to $C$ is oriented.
By Morse-Bott theory and the Thom isomorphism, there is a long exact sequence
$$ \cdots \to H_{S^1}^{*- \lambda}(C;\Q) \to H_{S^1}^{*}(M^+;\Q) 
\to H_{S^1}^{*} (M^-;\Q) \to  H_{S^1}^{*+1-\lambda}(C;\Q) \to  \cdots.$$
Moreover, the  composition  of the map  
$ H_{S^1}^{*-\lambda}(C;\Q) \to H_{S^1}^{*}(M^+;\Q)$
with the restriction map from  
$H_{S^1}^{*}(M^+;\Q) \to H_{S^1}^{*}(C;\Q)$ is precisely
multiplication by the Euler class of the negative normal bundle.
Since the circle fixes  $C$ but otherwise acts freely on
the negative normal bundle, the Euler class is not a zero divisor by 
Lemma~\ref{Atiyah-Bott}. Therefore, this composition is injective.
Thus the  map  
$ H_{S^1}^{*-\lambda}(C;\Q) \to H_{S^1}^{*}(M^+;\Q)$
is also injective.

In contrast, the function $\phi^2$ is not Morse-Bott at
the critical set  $C := \phi\inv(0)$. Nevertheless, 
because  $C$ is the minimum of $\phi^2$, the spaces $M^+$ and 
$C$ are homotopic, and $M^-$ is empty.  So the lemma follows immediately.
\end{pf}

Lemma~\ref{KM-circle} leads immediately to the following
theorems of Kirwan.

\begin{RefTheorem}
\labell{cor:injectivity-easy}
{\em (Kirwan)} Let $S^1$ act on a compact  symplectic manifold
$M$ with a moment map and let  $F$ denote
the fixed point set.
The natural restriction $H_{S^1}^*(M;\Q) \to  H_{S^1}^*(F;\Q)$ 
is  injective.
\end{RefTheorem}

This 
is a corollary of Lemma~\ref{KM-circle}
when we take the moment map as our Morse function.
It is proved by  induction on the proposition
that  the natural restriction 
$H_{S^1}^*(\phi\inv(-\infty,a);\Q) \to  
H_{S^1}^*(F~\cap~\phi\inv(-\infty,a);\Q)$ 
is an injection for a regular value  $a$ with $n$ critical values
below it.

\begin{RefTheorem}
\labell{cor:surjectivity-easy}
{\em (Kirwan)} Let $S^1$ act on a compact  symplectic manifold
$M$ with a moment map $\phi :M \to \R$.  
The Kirwan map $\kappa: H^*_{S^1}(M;\Q) \to H^*(M_\red;\Q)$
is surjective.
\end{RefTheorem}

This is also a corollary of 
Lemma~\ref{KM-circle} by an argument analogous to
Theorem~\ref{cor:injectivity-easy} when we take 
the function $\phi^2$ as our Morse function
and note that $0$ is its minimum value.

We can now state and prove our main theorem
for this special case.

\begin{Theorem}
Let $S^1$ act on a compact symplectic manifold $M$ with 
moment map $\phi: M \to R$.  Assume that $0$ is a regular value of $\phi$.
Let $F$ denote the set of fixed points; write $M_- := \phi^{-1}(-\infty,0)$,
and $M_+ := \phi^{-1}(0,\infty)$.
Define
$$K_\pm := \{ \alpha \in H_{S^1}^*(M;\Q) \mid \alpha|_{F \cap M_\pm} = 0 \},
\mbox{ and }$$
$$K := K_+ \oplus K_-.$$
Then there is a short exact sequence:
$$ 0 \to K \to H_{S^1}^*(M;\Q) \stackrel{\kappa}{\to}
 H^*(M_\red;\Q) \to 0,$$
where $\kappa: H^*(M;\Q) \to H^*(\Mred;\Q)$ is  the Kirwan map.
\end{Theorem}

\begin{Remark}
More generally, let a circle act on a manifold $M$.
A {\bf formal moment map} is a Morse-Bott function 
$\Phi: M \to \R$ such that
the critical points of $\Phi$ correspond exactly
to the fixed points.
(See  Ginzburg-Guillemin-Karshon \cite{GGK}.)
Then, as long as $M$ is compact and $0$ is a regular value of $\Phi$, 
the  theorem above is also true for any formal moment map; this
follows easily from a quick examination of the proof.
\end{Remark}

\begin{Remark}
Let $i_\pm : M_\pm \to M$ denote the natural inclusions.  By
Theorem~\ref{cor:injectivity-easy}, $K_\pm = {\rm ker}~i_\pm^*$.
In this case, 
a simple proof of our theorem can be given by 
an application of the Mayer-Vietoris sequence to the triple $(M,M_+, M_-).$
The proof given below is one which
generalizes to actions of higher-rank tori.
\end{Remark}

\begin{pf}
The injectivity of the inclusion $  K \to H_{S^1}^*(M;\Q)$ is obvious,
and the surjectivity of the Kirwan map 
$\kappa: H_{S^1}^*(M;\Q) \to H^*(\Mred;\Q)$ is 
Theorem~\ref{cor:surjectivity-easy}.

First, we show that  $K \subseteq \ker(\kappa)$.
By definition, for every  $\alpha \in K_+$, 
$\alpha|_{F \cap M_+} = 0$.
As in the proof of Theorem~\ref{cor:injectivity-easy},
applying Lemma~\ref{KM-circle} to the function $\phi$ and using induction, 
we see that $\alpha|_{\phi\inv(-\infty,\epsilon)} = 0$
for small $\epsilon > 0$.  In particular
$\alpha|_{\phi\inv(0)} = 0$.  Therefore, $K_+ \subseteq \ker(\kappa).$
A similar argument  using the function $ -\phi$ shows that
$K_- \subseteq \ker(\kappa).$  Therefore,
$K = K_+ \oplus K_- \subseteq \ker(\kappa)$.

It remains to show that $\ker(\kappa) \subseteq K$.
Order the set of connected components of
the fixed point set  $\{F_i\}_{i=1}^N$  of $M$  so that 
if $i < j$ then $\phi(F_i)^2 \leq \phi(F_j)^2$.
We wish to show that if  $\alpha \in H_{S^1}^*(M;\Q) $ satisfies
$\alpha|_{\phi\inv(0)} = 0$, then $\alpha$ is in $K$.
By Theorem~\ref{cor:injectivity} it suffices to
find $\alpha' \in K$ such that $\alpha|_{F_i} = \alpha'|_{F_i}$ for all $i$.
By induction, it suffices to prove that given  $p > 0$ and any
$\alpha \in H_{S^1}^*(M;\Q)$ which vanishes on $\phi\inv(0)$
and on $F_i$ for all $i < p$, there
exits $\alpha' \in K$ such that
$\alpha'|_{\phi\inv(0)} = \alpha|_{\phi\inv(0)}$
and $\alpha'|_{F_i} = \alpha|_{F_i}$ for all $i \leq p$.

Assume that we are given such an $\alpha$.
Applying  Lemma~\ref{KM-circle} to the function  $\phi^2$,
we see that $\alpha|_{F_p}$ is a multiple of the Euler class $e$ of the
negative normal bundle of $F_p$ for the function $\phi^2$.   
By symmetry, we may assume that $\phi(F_p) > 0$.
Then $F_p$ is also a critical set for the 
function $\phi$,
and  the Euler class of the negative normal bundle of
$F_p$ for $\phi$  is also $e$.  
Applying induction and Lemma~\ref{KM-circle} to the function $\phi$,
we see that there exists  
$\alpha' \in H_{S^1}^*(M;\Q)$ such that
$\alpha'|_{F_p} = \alpha|_{F_p}$,
and such that the restriction of $\alpha'$ to 
$\phi \inv(-\infty,f(F_p)-\epsilon)$ is trivial
for all $\epsilon > 0$. 
In particular, $\alpha'|_{F_i} = 0$ for all
$i$ such that $\phi(F_i) \leq 0$, and thus $\alpha' \in K_- \subset
K$.
Similarly, $\alpha'|_{F_i} = 0$ for all $i < p$.
The result follows.
\end{pf}

\section{Torus actions and Morse-Kirwan theory}

In this section, we review Morse-Kirwan theory 
and Kirwan's  application  of it to the cohomology ring  of
symplectic manifolds  with Hamiltonian torus actions.
The basic idea is that in order to generalize the earlier
result,  we do not really require our function to
be a non-degenerate Morse-Bott function. 
All we need is that our function
behave like a Morse function on the level of cohomology,
and that Lemma~\ref{Atiyah-Bott} holds for the normal
bundles of the critical sets.  Kirwan
shows that although functions arising from
the square of the moment map are not non-degenerate,
they do satisfy both these properties.

We have already encountered something of this kind in the
first section.
Let the circle $S^1$  act on a compact symplectic manifold
$M$ with moment map $\phi: M \to \R$.
As mentioned in the proof of Lemma~\ref{KM-circle},
the  critical set $\phi\inv(0)$ of function $\phi^2$ is
degenerate.  However, as $0$ is a minimum, the second
Morse lemma still holds at this critical level, so that
this kind of degeneracy is inconsequential.

Kirwan  developed an extension of Morse theory 
which applies to functions which  have well behaved degeneracies --
morally, they look like the product of  a minimum and a
non-degenerate Morse-Bott function.
In particular, she shows that for a critical set $C$ of such
a function $f: M \to \R$ 
there is a long exact sequence
$$ \cdots \to H_{S^1}^{*- \lambda}(C;\Q) \to H_{S^1}^{*}(M^+;\Q) 
\to H_{S^1}^{*} (M^-;\Q) \to  H_{S^1}^{*+1-\lambda}(C;\Q) \to  \cdots.$$
and that  the  composition  of the map  
$ H_{S^1}^{*-\lambda}(C;\Q) \to H_{S^1}^{*}(M^+;\Q)$
with the restriction map from  
$H_{S^1}^{*}(M^+;\Q) \to H_{S^1}^{*}(C;\Q)$ is precisely
multiplication by the Euler class of the negative normal bundle to $C$.
Here, $M^+ := f^{-1}(f(C) + \epsilon)$, and $M^- := f^{-1}(f(C) - \epsilon)$.
Kirwan then showed that this extension of Morse theory applies to
the square of the moment map for a torus action on a manifold.
Finally, she used this theory and  Lemma~\ref{Atiyah-Bott}
to argue exactly as in the proof of Lemma~\ref{KM-circle},
and thus prove the following Lemma:

\begin{Lemma} {\em (Kirwan)}
\labell{KM-torus}
Let a torus $T$ act on a symplectic manifold $M$ with a proper
moment map $\phi : M \to \ft^*$.
Choose a fixed inner product on the dual of the Lie algebra $\ft^*$.
Let a function be defined $f: M \to \R$ in either of the two following ways
\begin{enumerate}
\item
Given $a \in \ft^*$, define $f : M \to \R$ by 
$f(x) = \left< \phi(x) - a, \phi(x) - a \right>$.
\item
Given any $\xi \in \ft$,  define $f:= \phi^\xi$.
\end{enumerate}

Let $C$ be a critical set of index $\lambda$ for $f$ and assume that
$C$ is the only critical set in 
$f^{-1}( f(C) - \epsilon, f(C) + \epsilon)$ for some $\epsilon > 0$. 
Define $M^+ := f\inv(-\infty,f(C) + \epsilon)$, and 
$M^- := f\inv(-\infty,f(C) - \epsilon)$.
Then there exists a short exact sequence in equivariant cohomology
$$ 0 \to H_T^{*-\lambda}(C;\Q) \to 
H_T^{*} (M^+;\Q) \to H_T^{*}(M^-;\Q) \to 0.$$
Moreover, the composition of the injection from 
$H_T^{*-\lambda}(C;\Q) \to H_T^{*}(M;\Q)$
with  the restriction map from $H_T^{* }(M;\Q) \to H_T^{*}(C;\Q)$
is the cup product with the equivariant Euler class
of the negative normal bundle.
\end{Lemma}

The following theorems are direct consequences of this Lemma:

\begin{RefTheorem}
\labell{cor:injectivity}
{\em (Kirwan)} Let $T$ act on a symplectic manifold
$M$ with a moment map $\phi$.
Assume that the set of connected components of $F$, the fixed point set,
is finite,
and that there exists a
generic $\xi \in
\ft$  such that $\phi^\xi$ is proper and bounded below.
Then the natural restriction $H_{T}^*(M;\Q) \to  
H_{T}^*(F;\Q)$ is an injection.
\end{RefTheorem}

This follows by applying  Lemma~\ref{KM-torus} to 
the function $ \phi^\eta$ for a generic $\eta \in \ft$	
close to $\xi$, as in Theorem \ref{cor:injectivity-easy}.

Note that the restriction map may not be injective 
if we drop the assumption on
the image of $M$; for example, consider the cotangent bundle to $S^1$.

\begin{RefTheorem}
\labell{cor:surjectivity}
{\em (Kirwan)} Let $T$ act on a  symplectic manifold
$M$  with a proper moment map, and assume that the set of
connected  components of the fixed point set is finite.
If $0$ is a regular value of $\phi$, then
the Kirwan map $\kappa: H^*_{T}(M;\Q) \to H^*(M_\red;\Q)$
is a surjection.
\end{RefTheorem}

This follows by applying Lemma~\ref{KM-torus} to 
the function $||\phi||^2$, as in Theorem \ref{cor:surjectivity-easy}.

\section{The main theorem}

We now apply the ideas in the previous section to the 
case we are interested in: computing the cohomology of reduced spaces.

\begin{Theorem}\labell{main-theorem}
Let a torus $T$ act on a compact symplectic manifold $M$ with a
moment map $\phi: M \to \t^*$.  Assume $0$ is a regular value of $\phi$.
Let $F$ denote the set of fixed points.
For all $\xi \in \ft$, define
$$M_\xi := \{ m \in M \mid  \left< \phi(m), \xi \right> \leq 0 \} ,$$
$$K_\xi := \{ \alpha \in H^*(M;\Q) \mid \alpha|_{F \cap M_\xi} = 0\  \},
\hbox{ and}$$
$$K := \sum_{\xi \in \ft} K_\xi.$$
Then there is a short exact sequence:
$$ 0 \to K \to H_T^*(M;\Q) \stackrel{\kappa}{\to}
 H^*(M_\red;\Q) \to 0,$$
where 
$\kappa: H_T^*(M;\Q) \to H^*(\Mred;\Q)$ is  the Kirwan map.
\end{Theorem}
\begin{Remark}
As before, if we denote by $i_\xi$ the inclusion map $i_\xi : M_\xi \to M$,
then $K_\xi = {\rm ker~}i^*_\xi$.
\end{Remark}
\begin{pf}
The injectivity of the inclusion $K \to H_T^*(M;\Q)$ is obvious.
The surjectivity of the Kirwan map $ \kappa : H_T^*(M;\Q) \to H^*(M_\red;\Q)$ 
is Theorem~\ref{cor:surjectivity}.
It remains to show that $K = \ker( \kappa)$.

First we show that  $K \subseteq \ker(\kappa)$.
Choose $\xi \in \ft$ and $\alpha \in K_\xi$. 
By applying Lemma~\ref{KM-torus} to the function $f := \phi^\xi$ 
inductively, as in the proof of Theorem~\ref{cor:injectivity},
we get $\alpha|_{M_\xi} = 0.$ 
In particular, $\alpha|_{\phi\inv(0)} = 0$.
Since there are only a finite number of distinct $K_\xi$'s,
we are done. 

It remains to  show that $\ker(\kappa) \subset K$.
Order the critical sets $\{C_i\}_{i=1}^N$ of $||\phi||^2$ so that 
 $i < j$ exactly if $||\phi(C_i)||^2 < ||\phi(C_j)||^2$.
(We may assume
for simplicity that no two critical sets assume the
same value.)
We wish to show that if   $\alpha \in H_T^*(M) $ satisfies
$\alpha|_{\phi\inv(0)} = 0$, then $\alpha$ is in $K$.
By Theorem~\ref{cor:injectivity} it suffices to
find $\alpha' \in K$ such that $\alpha|_F = \alpha'|_F$.
By induction, it suffices to prove that given any $p > 1$ and any
$\alpha \in H_T^*(M)$ which vanishes on $C_i$ for all $i < p$, there
exists $\alpha' \in K$ such that
$\alpha'|_{C_i} = \alpha|_{C_i}$ for all $i \leq p$.
Here, we use the fact that every fixed set is a critical set,
and that the first critical set is $C_1 = \phi\inv(0)$.

Assume that we are given such an $\alpha$.
Applying Lemma~\ref{KM-torus}  to the function $||\phi||^2$,
we see that $\alpha|_{C_p}$ is a multiple of the Euler class $e$ of the
negative normal bundle of $C_p$  for the function $||\phi||^2 $. 
The point $x \in M$ is critical for $||\phi-a||^2$ exactly
if the vector $(\phi(x)-a)_M$ is zero,
where for any $b \in \ft^*$, $b_M$ denotes the
vector field on $M$ associated to the one parameter
subgroup generated 
by the element in $\ft$ associated to $b$ by the invariant inner product.
Therefore, for any  $\lambda \in \R^+$, 
$C_p$ is  a critical set for the 
function $||\phi + \lambda\phi(C_p)||^2$.
Moreover, the
Euler class of the negative normal bundle on $C_p$ for this new function
is still $e$.

Since $M$ is compact, 
for sufficiently large  $\lambda$,
$||\phi(F) + \lambda\phi(C_p)||^2
< ||\phi(C_p) + \lambda\phi(C_p)||^2$
for  all $F \in \cF$ whose inner product with $C_p$ is negative.
Choose such a $\lambda$ and apply Lemma~\ref{KM-torus} 
to the function $f: = ||\phi~+~\lambda\phi(C_p)||^2$.
Since $\alpha|_{C_p}$ is a multiple of $e \in H_T^*(C_p)$,
there exists $\alpha' \in H^*(M)$ such that
$\alpha'|_{C_p} = \alpha|_{C_p}$ 
and the restriction $\alpha|_{C_i}$ is trivial for all
$i$ such that $f(C_i) < f(C_p)$.
In particular, $\alpha'|_F = 0$ for all $C_i$ whose  
inner product with $\phi(C_p)$ is negative. Hence $\alpha'$
lies in $K$.
Finally, for all $i < p$,
$||\phi(C_i) + \lambda\phi(C_p)||^2 < ||\phi(C_p) + \lambda\phi(C_p)||^2,$
and hence $\alpha'|_{C_i} = 0$.
\end{pf}

\begin{Remark} ({\bf Chern classes})
Theorem \ref{main-theorem}  allows the cohomology ring of the
reduced space $\Mred$ to be computed 
in terms of fixed point data on $M$, namely the  cohomology ring
of the fixed point set and the restriction of each equivariant
cohomology class on $M$ to this set.
The Chern classes of the reduced space $\Mred$ can also be computed in 
terms of fixed point data, namely, the restriction of each equivariant
Chern class to the fixed point set.
This is because $TM|_{\phi^{-1}(0)} = \pi^*T\Mred
\oplus \ft_{\C}^*$, where $\pi: \phi^{-1}(0) \to \phi^{-1}(0)/T = \Mred$
is the projection, and $\ft_{\C}^*$ is a trivial bundle.  Thus, in terms
of the Kirwan map $\kappa$, $c_i(\Mred) = \kappa (c_i(M))$.  On the other
hand, by injectivity, $c_i(M)$ is determined by its image in $H_T^*(F)$.
\end{Remark}

\begin{Remark}({\bf Computability})
The description of the kernel of $\kappa$ given in Theorem~\ref{main-theorem}
is 
algorithmically computable in the following sense.  The cohomology ring
$H_T^*(M)$ is a finitely generated module over $H_T^*({\rm pt})$.  Suppose
we are given a finite basis for this module.  Then by solving a finite
set of 
linear equations, we can produce a finite basis for the submodule
$K_\xi$, for any $\xi$.  Since only a finite number of elements $\xi \in \ft^*$
give rise to distinct submodules $K_\xi$, we may produce a finite basis
for $K$.
\end{Remark}

\begin{Theorem}
\labell{proper-theorem}
Let a torus $T$ act on a  symplectic manifold $M$ with a
moment map $\phi: M \to \R$.  Assume $0$ is a regular value.
Assume that   the set  of connected components of the
set of fixed points is finite and that there exists 
$\eta \in \ft$ such that $\phi^\eta$ is proper and bounded below.
For all $\xi \in \ft$, define
$$M_\xi := \{ m \in M \mid  \left< \phi(m), \xi \right> \leq 0 \} ,$$
$$K_\xi := \{ \alpha \in H^*(M;\Q) \mid \alpha|_{M_\xi} = 0\  \} , and$$
$$K := \sum_{\xi \in \ft} K_\xi.$$
Then there is a short exact sequence:
$$ 0 \to K \to H^*(M;\Q) \stackrel{\kappa}{\to}
 H^*(M_\red;\Q) \to 0,$$
where 
$\kappa: H_T^*(M;\Q) \to H^*(\Mred;\Q)$ is  the Kirwan map.
\end{Theorem}

The argument in the non-compact case is nearly identical to
the proof in the compact case; except
that since $\phi^\xi$ may not be proper, there may exist
cohomology classes which vanish on  $M_\xi \cap F$ but
not on $M_\xi$  itself.

\section{Extending to the integers}

So far in this paper we have restricted our attention to rational
cohomology.  In fact, many of the results we have proved have their
analogs in integral cohomology; but there are a number of subtleties
which occur.

In order to use any of the methods of this paper, we
must verify that an analog of Lemma~\ref{Atiyah-Bott} holds over the
integers.  It is 
easy to see that
the naive generalization
of Lemma~\ref{Atiyah-Bott} is false.
Consider, for example, the trivial plane bundle $\C \oplus \C$ 
over $\RP^3$,
and let $S^{1}$ act on the plane bundle by $\lambda \cdot (y,z)  =
(\lambda y, \lambda^2 z)$.  Since $S^1$  acts trivially on
${\RP}^3$, the  equivariant cohomology of
$\RP^3$ is the tensor product of the cohomology of $\RP^3$ and
the cohomology of $BS^1 = \CP^\infty$.  The Euler class
of this bundle is $2x^2$, where $x$ denotes the generator of  the
cohomology of $\CP^\infty$.  This class clearly is a zero divisor, as
it annihilates the two dimensional cohomology of $\RP^3$.

In order to avoid such examples, we have to place some restrictions on
the action, or else on the cohomology of the fixed point set.  The local
result we obtain is the following.

\begin{Lemma}\labell{A-B;Z}
Let a torus $T$ act on 
an oriented vector bundle $E$ over a manifold $X$, fixing $X$.
Assume that there exists $S^1 \subset T$ which
does not fix any point in $E \ssminus X$.
Assume also that,  for each prime $p$, 
one of the following two conditions is satisfied:
\begin{enumerate}
\item The cohomology of $X$ has no $p$ torsion.
\item For every point  $m \in E \ssminus X$,
there exists a $\Z/p \subset  T$  which acts freely on $m$. 
\end{enumerate}
Choose an invariant metric on $E$, and let $D$ and $S$ denote the unit disk and
sphere bundle, respectively.
The  natural long exact sequence in relative equivariant cohomology
splits into a short exact sequence:
$$ 0 \to H^*(D,S;\Z) \to  H^*(D;\Z) \to H^*(S;\Z) \to 0. $$
\end{Lemma}

Again, this is equivalent to the statement that the Euler
class of $E$ is not a zero divisor.
To prove this, it suffices to show that each weight 
is not a zero divisor.
Therefore, it is enough to
note that for any $T$-representation on $\C$,
the corresponding weight is not divisible by $p$ if
there exists a $\Z/p \subset  T$  which acts freely on every point
except $0$. 

Checking the proof of Kirwan's injectivity theorem (Theorem
\ref {cor:injectivity}),
it is easy to see that it holds over the integers wherever this
version of the Atiyah-Bott lemma holds for every fixed point. 
Thus:

\begin{Proposition}
\labell{cor:injectivity-Z}
Let a torus $T$ act on a symplectic manifold
$M$ with a moment map $\phi$.
Assume that the set of connected components of $F$, the fixed point set,
is finite, and that there exists a
generic $\xi \in
\ft$  such that $\phi^\xi$ is proper and bounded below.
Suppose in addition that 
for every prime $p$, one of the following two conditions is satisfied:
\begin{enumerate}
\item  The integral cohomology of $F$ has no $p$-torsion, or:
\item For every point $m \in M$ which is not fixed by the $T-$action,
there exists a subgroup of $T$ congruent to $\Z/p$ which acts freely
on $m$.
\end{enumerate}
Then the natural restriction $H_{T}^*(M;\Z) \to 
H_{T}^*(F;\Z)$ is an injection.
\end{Proposition}

In particular, injectivity holds
if the fixed point set has no torsion, or if the
action is quasi-free. 

Integer analogues of Kirwan's surjectivity theorem (Theorem
\ref{cor:surjectivity}),
and our main result (Theorem \ref{main-theorem}) also hold
under a  similar set of  restrictions in the case of an $S^1$ action. 
Once again, the proof amounts to noting that
the proofs given earlier go through  without change as long as
the Atiyah-Bott lemma holds for the negative normal bundle
at each critical point. In this case,  
each  critical point is either fixed, in which case we apply  
Lemma \ref{A-B;Z} above, or the minimum, which works automatically.

\begin{Proposition}
\labell{cor:surjectivity-Z}
Let $S^1$ act on a compact  symplectic manifold
$M$ with a moment map $\phi :M \to \R$.
Assume that 
for every prime $p$, one of the following two conditions is satisfied:
\begin{enumerate}
\item  The integral cohomology of $F$ has no $p$-torsion, or:
\item For every point $m \in M$ which is not fixed by the $T-$action,
there exists a subgroup of $T$ congruent to $\Z/p$ which acts freely
on $m$.
\end{enumerate}
Then the Kirwan map $\kappa: H^*_{S^1}(M;\Z) \to H^*(M_\red;\Z)$
is surjective.

\end{Proposition}

Likewise,

\begin{Proposition} \labell{S1-Z}
Let $S^1$ act on a compact symplectic manifold $M$ with 
moment map $\phi: M \to R$.  Assume that $0$ is a regular value.
Let $F$ denote the set of fixed points.
Assume that 
for every prime $p$, one of the following two conditions is satisfied:
\begin{enumerate}
\item  The integral cohomology of $F$ has no $p$-torsion, or:
\item For every point $m \in M$ which is not fixed by the $T-$action,
there exists a subgroup of $T$ congruent to $\Z/p$ which acts freely
on $m$.
\end{enumerate}

Define
$$K_+ := \{ \alpha \in H_{S^1}^*(M;\Z) \mid \alpha|_{F_+} = 0 \},
\mbox{    where }
F_+ := F \cap \phi\inv(0, \infty); \mbox{    }$$ 
$$K_- := \{ \alpha \in H_{S^1}^*(M;\Z) \mid \alpha|_{F_-} = 0 \},
\mbox{    where }
F_- := F \cap \phi\inv(-\infty, 0); \mbox{ and} $$
$$K := K_+ \oplus K_-.$$

Then there is a short exact sequence:
$$ 0 \to K \to H_{S^1}^*(M;\Z) \stackrel{\kappa}{\to}
 H^*(M_\red;\Z) \to 0,$$
where $\kappa: H^*(M;\Z) \to H^*(\Mred;\Z)$ is  the Kirwan map.
\end{Proposition}

\begin{Remark} In the case where $S^1$ does not act freely
on $\phi\inv(0)$, the cohomology group $H^*(M_{\rm red}, \Z)$ which
enters into both these results is 
the {\em orbifold} cohomology of $M_\red$, not the cohomology of the
underlying topological space.  These are of course identical over $\R$
but are, in general, different over $\Z$.
\end{Remark}

On the other hand, 
when a higher-rank torus $T$ acts on a manifold, 
Kirwan's surjectivity Theorem \ref{cor:surjectivity},
and our main result (Theorem \ref{main-theorem}) require a more subtle
application of the Atiyah-Bott Lemma, in that the critical points of the
Morse function $\phi^2$ include spaces more complicated than the fixed
points.  In order to generalize our proofs, we  need a version of
Lemma \ref{A-B;Z} which holds at all of these critical points.  There
are various conditions under which this can be seen to hold, but the
statements of these conditions can be cumbersome.
We therefore restrict our attention
to the quasi-free case, where the Euler class of the negative
normal bundle (at any critical manifold) is
not a zero divisor; thus the conditions of Lemma 6.1 are met at
every point of the fixed set.

\begin{Proposition}
\labell{cor:surjectivity-Z-torus}
{\em (Kirwan)} Let $T$ act quasi-freely on a  symplectic manifold
$M$  with a proper moment map, and assume that the set of
connected  components of the fixed point set is finite.

If $0$ is a regular value, then
the Kirwan map $\kappa: H^*_{T}(M;\Z) \to H^*(M_\red;\Z)$
is a surjection.
\end{Proposition}

Finally we have an integer version of our main result:

\begin{Proposition}
Let a torus $T$ act quasi-freely on a  symplectic manifold $M$ with a
moment map $\phi: M \to \R$.  Assume $0$ is a regular value.
Assume that   the set   of connected components of $F$, the
set of fixed points,  is finite and that there exists 
$\xi \in \ft$ such that $\phi^\xi$ is proper and bounded below.
For all $\xi \in \ft$, define
$K_\xi$ to be the set of $ \alpha \in H^*(M;\Z)$ such that
$\alpha$ vanish when restricted to the set $M_\xi= (\phi^\xi)^{-1}(0,\infty)$.
Define $K := \sum_{\xi \in \ft} K_\xi$.
Then there is a short exact sequence:
$$ 0 \to K \to H^*(M;\Z) \stackrel{\kappa}{\to}
 H^*(M_\red;\Z) \to 0,$$
where 
$\kappa: H_T^*(M;\Z) \to H^*(\Mred;\Z)$ is  the Kirwan map.
\end{Proposition}

\section{Torsion-free reduced spaces}

In a large class of examples, the integer cohomology
of the reduced spaces is torsion-free.  We begin with
the case of reduction by an $S^1$ action.

\begin{Theorem}\labell{thm:torsion-free}
Let the circle $S^1$ act on a compact symplectic manifold $M$ with
a moment map $\phi$.  Assume that the integer cohomology 
$H^*(F;\Z)$ of the fixed point set $F$ is torsion-free, and that $S^1$ acts
freely on $\phi\inv(0)$.  Then the integer cohomology ring
of the reduced space is also torsion-free. 
\end{Theorem}

\begin{pf}

Note first that because the fixed point set is torsion-free,
the conditions of Propositions 6.2 and 6.3 are met.
Assume that $p \in \Z$ is a prime and $\alpha
\in H^*_{S^1}(M;\Z)$ is a cohomology class and $p \alpha|_{\phi\inv(0)} =
0.$ 
We must show that  $\alpha|_{\phi\inv(0)} = 0$.
Order the set of connected components of
the fixed point set  $\{F_i\}_{i=1}^N$  of $M$  so that 
if $i < j$ then $\phi(F_i)^2 \leq \phi(F_j)^2$.
By Theorem~\ref{cor:injectivity} it suffices to
find $\alpha' \in H_{S^1}^*(M;\Z)$ such that 
$\alpha'|_{\phi\inv(0)} = 0 $ and $\alpha_{F_i} = \alpha'|_{F_i}$ for all $i$.
By induction, it suffices to prove that given an integer  $q > 0$ and any
$\alpha \in H_{S^1}^*(M;\Z)$ which vanishes on  
$F_i$ for all $i < q$ and such that
  $p \alpha$ vanishes on $\phi\inv(0)$, there
exists $\alpha' \in H_{S^1}^*(M;\Z)$  such that   
and $\alpha'|_{F_i} = \alpha|_{F_i}$ for all $i \leq q$
and  $\alpha'$ vanishes on $\phi\inv(0)$. 

Assume that we are given such an $\alpha$.
Applying  Lemma~\ref{KM-circle} to the function  
$\phi^2$ and the class $p \alpha$,
we see that $p \alpha|_{F_q}$ is a multiple of the Euler class $e$ of the
negative normal bundle of $F_q$ for the function $\phi^2$.   
We now have two cases two consider.

First, assume that $p \not|\ e$.
Then, by Lemma~\ref{Technical}, $\alpha|_{F_q}$ is also a multiple of $e$. 
By symmetry, we may assume that $\phi(F_q) > 0$.
Then $F_q$ is also a critical set for the 
function $\phi$,
and  the Euler class of the negative normal bundle of
$F_q$ for $\phi$  is still $e$.  
Applying induction and Lemma~\ref{KM-circle} to the function $\phi$,
we see that there exists  
$\alpha' \in H_{S^1}^*(M;\Q)$ such that
$\alpha'|_{F_q} = \alpha|_{F_q}$,
and such that the restriction of $\alpha'$ to 
$\phi \inv(-\infty,f(F_q)-\epsilon)$ is trivial
for all $\epsilon > 0$. 
In particular, 
$\alpha'|_{F_i} = 0$ for all $i < p$.
and $\alpha'$ vanishes on $\phi\inv(0)$.

More generally,  
let $N \subset M$ be the connected component of the set of
points which is fixed by $\Z/(p)$ which contains $F_q$.
Let $a$ be the downward Euler class for $F_q$ in $N$ using
the Morse-Bott function $\phi$,
an $b$ be the Euler class for the complimentary bundle, so
$e = a b$.
Since the action on $\phi\inv(0)$ is free, the induction
hypothesis assures that $\alpha|_{F_i} = 0$  for
all $i$ such that $\phi(F_i) < \phi(F_q)$ and $F_i \subset N$
Therefore $\alpha|_{F_q} = a \eta$, for some $\eta$.
Since $p \alpha|_{F_q} = e \xi = a b \xi$,
we can substitute to get $p a \eta = a b \xi$.
Since $a$ cannot be a zero-divisor, this implies
that $p \eta = b \xi$.  As before, this implies that
$\xi = p \xi'$, etc.
\end{pf}

\begin{Lemma} 
\label{Technical}
Let a compact torus  $T$ act
on an oriented real vector bundle $E$ over a manifold $X$.
Assume that a circle subgroup $S^1 \subset T$ acts on $E$
so that the fixed point set of $S^1$ and $\Z/(p) \subset S^1$
is precisely $X$ for $p \in \bN$.
Assume also that $H^*(X;Z)$ is torsion-free.
Let $e$ be the equivariant Euler class of $E$.
Then for any cohomology classes $\alpha$ and $\beta$ such that
$p \alpha = \xi e$, there  exists $\xi'$ such that $p \xi' = \xi$.
\end{Lemma}

We start by assuming that $X$ is fixed by $T$.

\begin{pf}
Identify the torus $T$ with $(S^1)^n$.
Order $n$-tuples as follows:
first by the sum of  all coordinates,
then by the first coordinate, the second coordinate, etc. 
Because the action fixes $X$, 
the cohomology ring $H_T^*(X;\Z)$ is simply the polynomial
ring in $n$ generators on $H^*(X;\Z)$.
Every term  of $e$ such that the sum of the coefficients is $\lambda$
is an integer, where $\lambda$ is the rank of $E$. 
Because the action  is $p$-free outside $X$,  
$p$ cannot divide all these integers.
Let $I$ be the largest $n$-tuple such that $p \not| \xi_I$,
and $J$ be the largest $n$-tuple such that $p \not| e_J$.
Note that be the above discussion $e_J$ is an integer.
Now consider the coefficients of the term $X^{I+J}$ in
the expression $p \alpha = e \xi$.
It is $$\alpha_{I+J} = \sum_{I' + J' = I + J} e_{I'} \xi_{J'}.$$
But if $I' + J' = I + J$, then either $I = I'$ and $J = J'$ or
$I' > I$ or $J' > J$.  Therefore, $e_I \xi_J = p \beta$,
for some $\beta$.
We can find integers $s$ and $t$ so that $1 = sp + te_I$.
Therefore, $\xi_J = sp \xi_J + t e_I \xi_J = p (s \xi_J + t \beta)$.
\end{pf}

 \begin{Remark}
Note that if $S^1$ does not act freely on $\phi\inv(0)$, then the
 the cohomology of the reduced space must have torsion.
 \end{Remark}

By examining the proof of Theorem \ref{thm:torsion-free}, we see
what difficulties will arise in the case of reduction by the action
of a torus of rank $r > 1$:  we
require surjectivity of the Kirwan map, which, in turn,
requires a condition on either the stabilizers of points or on the
cohomology of critical points of the square of the moment map.  These
are themselves reduced spaces, whose cohomology being torsion-free
will depend on a version of Theorem \ref{thm:torsion-free} for reduction
by tori of rank $r-1$.  Thus the reduced spaces will be torsion-free
under a complicated set of conditions, defined inductively.  These 
conditions can easily be verified in examples but are not so easy
to state in full generality; we content ourselves with the following
statement.

\begin{Theorem}\labell{thm:torsion-free-qf}
Let a torus $T$ act quasi-freely
on a symplectic manifold $M$ with moment map $\phi$.  Suppose $0$ is a regular value
of the moment map, and that the integer cohomology $H^*(F;\Z)$ of the
fixed point set $F$ is torsion-free.
Then the cohomology $H^*(\Mred;\Z)$ of the reduced space $\Mred$
is also torsion-free.
\end{Theorem}
 
\section{Example:  Smooth Toric Varieties}

In this section, 
to demonstrate our theorem,
we show how it provides a simple method to obtain the cohomology 
rings of smooth compact projective toric varieties
(see e.g. \cite{Dan}).  These compact symplectic 
manifolds can all be obtained as symplectic quotients of
Euclidean spaces.

The moment map $\psi :\C^N \to {\R^N}^*$ 
for the natural action of $(S^1)^N$ on $\C^N$ is 
$\psi(z_1,\ldots,z_N) = (|z_1|^2, \ldots, |z_N|^2)$. 
The image of this moment map is given by
$\psi(\C^N)= \{\xi \in {\R^N}^* \mid \xi_i \geq 0 \mbox{ for all } i \}$.

Consider a subgroup $G \subset (S^1)^N$.
This inclusion gives rise to a short exact sequence of Lie algebras 
$$0 \to \fg \stackrel{i}{\to} \R^N  \stackrel{\pi}{\to} \ft \to 0$$
and its dual sequence
$$0 \to \ft^* \stackrel{\pi^*}{\to} {\R^N}^* \stackrel{i^*}{\to} \fg^* \to 0.$$

For any $\eta \in \fg^*$, the function $\phi := i^* \circ \psi - \eta$
is a  moment map  for the induced action of $G$ on $\C^N$. 
Assume that $\phi$ is proper, and that $0$ is a regular value of $\phi$. 
Let $M$ be the reduction of  $\C^N$ by $G$ at $0$. 
The torus $T := (S^1)^N/G$ acts symplectically on $M$.
The image of moment map for the $T$ action
can be identified with the polytope 
$\Delta := \{\xi \in {\R^N}^* \mid \xi_j \geq 0 \mbox{ for all } j 
\mbox{ and } i^*(\xi) - \eta = 0 \}.$

There is a (possibly empty) facet of $\Delta$ for each $i \in
(1,\ldots,N)$, given by 
$\xi  \in {\R^N}^*$ such that $ \xi_j \geq 0$  for all  $j$ ,
$i^*(\xi) - \eta  = 0$ and $\xi_i = 0$.
For $I \subset (1,\ldots,N)$, we say that the $I$ facets
{\em intersect} if the intersection of these sets is not
empty.

Since $  \dim T =\frac{1}{2} \dim M$, $M$ is a compact 
smooth toric variety.

\begin{Theorem}
Let the ideals $\cI$ and $\cJ \subset \Q[x_1,\ldots,x_N]$
be defined as follows:
 $\cJ = \{ \sum \alpha_i x_i \mid \alpha \in \pi^*(\ft^*) \}$,
and  
$\cI$ is the ideal
 generated by $\prod_{i \in I} x_i$ for all $I \subset (1,\ldots,N)$
such that
the $I$ facets do not  intersect in
the polytope $\Delta$.
The cohomology ring $H^*(M;\Q)$ is given by
$\Q[x_1,\ldots,x_N]/(\cI,\cJ)$.
\end{Theorem}

\begin{pf}
The $G$-equivariant cohomology of $\C^N$ can be expressed
as $$H_G^*(\C^N;\Q) = \Q[x_1,\ldots,x_N]/\cJ.$$

From the construction of $\phi$, it is clear that since $\phi$
is proper, there exists $\zeta$ such that $\phi^\zeta$ is proper
and bounded below.  Therefore, we  may apply
Theorem~\ref{proper-theorem}.

Consider a generic $\xi \in \fg$.
The only critical point of $\phi^\xi$ is $0$,
and the only critical value is $-\left<\eta,\xi\right>$.

If $\left< \eta, \xi \right> > 0$, then $0 \in M_\xi$.
Since every cohomology class which vanishes when restricted
to $0$ is trivial, $K_\xi = 0$.

Now consider the case $\left< \eta, \xi \right> < 0$.
The subbundle of the tangent space
at $0$ consisting of points with all the coordinates
zero except the $i$'th is in the negative normal bundle  for
$\phi^\xi$  exactly if 
 $\left< \beta_i, \xi \right> < 0 $,
where $\beta_i = \pi^*(e_i)$, and the $e_i$'s are the standard basis of
${\R^N}^*$.  
Therefore, $K_\xi$ is generated by $\prod_{i \in I} x_i$,
where
$I := \{ i \in (1,\ldots,N) \mid \left< \beta_i, \xi \right> < 0\}$.
(We know that its restriction to $0$ must be of this form, and
every class is determined by its restriction to $0$.)

Putting this together, it is clear that $K$ is generated
by the products $\prod_{i \in I} x_i$,
over all sets $I$ such that  
there exists $\xi$ with $\left< \eta, \xi \right> < 0$ and 
$\left< \beta_j, \xi \right> > 0$ for all $j \not\in I$. 

Given $I \subset (1,\ldots,N)$ 
it is clear that $\eta$ is not in $\sum_{j \not\in I} \R^+ \beta_j$
exactly if there exists $\xi$ such that
$\left< \eta, \xi \right> < 0$ and 
$\left< \beta_j, \xi \right> > 0$ for all $j \not\in I$. 
Thus, $K$ is generated by $\prod_{i \in I} x_i$ for all sets
$I$ such that $\eta$ is not in 
$\sum_{j \not\in I} \R^+ \beta_j$.
Finally, $\eta$ is in $\sum_{j \not\in I} \R^+ \beta_j$
exactly if the facets corresponding to the elements of $I$ 
intersect in $\Delta$.
\end{pf}

\end{document}